\documentclass[brochure,12pt]{bourbaki}
\usepackage[matrix,arrow]{xy}
\usepackage{amssymb,amsfonts,amsmath,footnote}
\usepackage[francais]{babel}
\usepackage{graphicx}
\usepackage[utf8]{inputenc}
\newcommand{\C}{{\mathbb C}}
\newcommand{\HH}{{\mathbb H}}
\newcommand{\HP}{{\mathbb H\mathbb P}}

\newcommand{\R}{{\mathbb R}}
\newcommand{\Z}{{\mathbb Z}}
\newcommand{\cA}{{\mathcal A}}

\newcommand{\cF}{{\mathcal F}}
\newcommand{\cG}{{\mathcal G}}
\newcommand{\cL}{{\mathcal L}}

\newcommand{\cP}{{\mathcal P}}

\newcommand{\cT}{{\mathcal T}}

\newcommand{\isom}{{\rm Isom}}
\newcommand{\Adm}{{\rm Adm}}
\newcommand{\adm}{{\rm adm}}
\newcommand{\Hom}{{\rm Hom}}

\date{juin 2015}
\bbkannee{67\`eme ann\'ee, 2014-2015}
\bbknumero{1103}
\title{Variétés lorentziennes plates\\ vues comme limites de variétés anti-de Sitter}
\subtitle{d'apr\`es Danciger, Guéritaud et Kassel}
\author{Jean-Marc SCHLENKER}
\address{Universit\'e du Luxembourg\\ 
Campus Kirchberg\\ 
Unit\'e de Recherche en Math\'ematiques\\
6, rue Richard Coudenhove-Kalergi\\
L-1359 Luxembourg} 
\email{jean-marc.schlenker@uni.lu}

\begin{document}
\maketitle


\noindent{\bf INTRODUCTION}

\bigskip


Les travaux de Danciger, Guéritaud et Kassel 
\cite{danciger-gueritaud-kassel:geometry,danciger-gueritaud-kassel:margulis} 
développent de manière féconde
les liens entre trois directions de recherche~:
\begin{itemize}
\item les variétés lorentziennes plates complètes de dimension $3$, en relation avec la conjecture
de Margulis qui les concerne,
\item les variétés anti-de Sitter complètes de dimension $3$,
\item les applications contractantes entre surfaces hyperboliques, ou plus généralement les 
applications contractantes équivariantes du plan hyperbolique dans lui-même, et leurs versions
infinitésimales, les champs de vecteurs contractants sur les surfaces hyperboliques.
\end{itemize}
Nous allons présenter rapidement quelques résultats obtenus par Danciger, Guéritaud
et Kassel dans ces différents domaines, en 
commençant par les énoncés sur les surfaces avant de passer à la dimension 3. Nous
verrons alors comment les résultats concernant les surfaces ont des interprétations
naturelles en termes de variété de dimension 3. 

Il n'est pas question ici de donner une présentation exhaustive des preuves, d'autant
que les références \cite{danciger-gueritaud-kassel:geometry,danciger-gueritaud-kassel:margulis}
sont remarquablement bien écrites. On va tenter de présenter de manière assez synthétique
les principaux résultats, et de mettre en évidence l'articulation entre eux, ainsi que 
leur contexte. On se concentrera en particulier sur 
\begin{itemize}
\item la notion de {\em géométrie
transitionnelle}, qui offre un pont entre les géométries hyper\-bolique, Minkowski et
anti-de Sitter, que les auteurs utilisent pour donner des preuves simples de résultats
importants sur les variétés lorentziennes plates à partir de résultats correspondants
en géométrie anti-de Sitter, dans le $\S$\ref{sc:transition},
\item la relation entre feuilletages par des géodésiques de type temps et
difféomorphismes contractants (resp. champs de vecteurs contractants), 
dans le $\S$\ref{sc:geodesiques},
\item la relation entre déformations en bandelettes (resp. déformations infinitésimales 
en bandelettes) de surfaces hyperboliques et domaines fondamentaux bordés par des
plans croches, dans le $\S$\ref{sc:defos-croches}.
\end{itemize}

\subsubsection*{Notations} On notera $S$ une surface, homéomorphe à l'intérieur d'une
surface compacte à bord, munie d'une métrique hyperbolique convexe co-compacte (c'est-à-dire
dont tous les bouts sont d'aire infinie). On appellera $\cT_S$ l'espace de Teichm\"uller
de $S$, vu comme espace des métriques hyperboliques convexes co-compactes sur $S$ considérées
à isotopie près.

\section{Le complexe des arcs et les déformations par bandelettes}

On s'intéresse ici à l'espace des métriques hyperboliques complètes sur une surface~$S$ 
homéomorphe à l'intérieur d'une surface compacte orientée à bord. Plus spécifiquement, on
peut se poser la question suivante~: comment peut-on passer de manière {\og géométrique\fg}
d'une métrique hyperbolique $h$ fixée à une autre métrique hyperbolique $h'$ ? 

Si on se place du point de vue des structures complexes (ou conformes) sur $S$, une réponse est donnée
par la théorie {\og de Teichm\"uller\fg} des applications quasi-conformes optimales entre 
surfaces de Riemann. 

Du point de vue des surfaces hyperboliques, lorsque $S$ est fermée (compacte sans bord)
une autre réponse possible est donnée par les coordonnées de Fenchel-Nielsen (voir 
par exemple \cite{fenchel-nielsen}) qui ont l'inconvénient de dépendre du choix d'une décomposition
en pantalons. Si on veut éviter un tel choix, on peut utiliser le théorème de tremblement de terre de 
Thurston (démontré dans \cite{kerckhoff}) qui affirme qu'il existe un unique tremblement
de terre qui relie $h$ à $h'$ --- les tremblements de terre sont des extensions aux 
laminations mesurées des twists de Dehn fractionnels le long de multicourbes pondérées. Il est d'ailleurs
intéressant de constater que la preuve la plus simple du théorème de tremblement de terre
de Thurston, obtenue par Mess \cite{mess,mess-notes}, utilisait déjà la géométrie anti-de Sitter.

Une autre alternative possible, pour les surfaces fermées, est fournie par l'application 
de grafting, qui associe à une métrique hyperbolique $h$ et à une lamination mesurée $l$ une
autre métrique hyperbolique $gr(m,l)\in \cT$ sur $S$. Lorsque $l$ est une courbe
fermée~$c$ munie d'un poids $w>0$, la structure conforme sous-jacente à $gr(m,l)$ 
est obtenue en coupant $m$ le long de $c$ et en y introduisant une bande plate de largeur $w$. 
\mbox{Dumas} et Wolf \cite{dumas-wolf} ont montré que, si $m$
est fixée, l'application $l\mapsto gr(m,l)$ est un homéomorphisme, fournissant ainsi une
autre paramétrisation de $\cT_S$.

Le théorème des tremblements de terre s'étend aux surfaces hyperboliques complètes \cite{mbh}, 
mais le résultat de Danciger, Guéritaud et Kassel permet de paramétrer spécifiquement les
déformations {\em expansives} par rapport à une représentation $\rho$ donnée. 
Ils définissent une forme particulière de déformation 
des surfaces hyperboliques qui, par construction, augmentent la longueur de toutes
les courbes fermées.

\subsection{Déformations par bandelette}

On considère une surface hyperbolique complète $S$, homéomorphe à l'intérieur
d'une surface compacte à bord $\bar S$. On supposera que tous les bouts de $S$ sont d'aire
infinie, et ne sont donc pas des cusps, $S$ est alors dite {\em convexe co-compacte}.

\begin{defi}
Un {\em arc géodésique} de $S$ est une géodésique complète dont chaque bout sort de tout 
compact de $S$. 
\end{defi}

Alternativement, on considérera la notion topologique correspondante, et on 
appellera arc sur $S$ une classe d'homotopie de courbe tracée sur $\bar S$, dont
les extrémités se trouvent sur le bord. Néanmoins un arc géodésique de $S$ n'est 
pas déterminé uniquement par l'arc (topologique) correspondant, puisqu'on peut
déplacer ses extrémités sur le bord à l'infini de $S$ sans changer d'arc topologique. 

\begin{defi}
On appellera {\em bandelette hyperbolique} la région du plan hyperbolique bordée par
deux géodésiques ultraparallèles, c'est-à-dire qui ne se rencontrent ni dans le
plan hyper\-bolique $\HH^2$, ni dans son bord idéal. La {\em taille} d'une 
bandelette hyperbolique est le segment géodésique qui joint ses deux composantes
de bord en étant orthogonal à chacune. La {\em largeur} d'une bandelette est
la longueur de sa taille.
\end{defi}

On notera qu'une bandelette hyperbolique est uniquement déterminée, à isométrie
près, par sa largeur. 

\begin{defi}
Une {\em déformation par bandelette} de $S$ le long de l'arc géodésique 
$\alpha$ est la surface hyperbolique complète obtenue en coupant $S$ le 
long de $\alpha$ et en y recollant une bandelette hyperbolique $b$, de manière
que les deux extrémités de la taille de $b$ soient identifiés au même point 
de $\alpha$.
\end{defi}

On note qu'une déformation par bandelette est uniquement déterminée si, en 
plus de~$S$ et de $\alpha$, on se donne la largeur $w_\alpha$ de la 
bandelette qu'on recolle, et le point $p_\alpha$ de~$\alpha$ où se trouveront les
extrémités de cette taille.

Si $\alpha_1,\cdots, \alpha_k$ sont des arcs géodésiques disjoints, on dira
qu'une déformation par bandelettes de $S$ le long de $\alpha_1,\cdots,\alpha_k$
est une surface hyperbolique complète obtenue en faisant successivement des
déformations par bandelette le long de $\alpha_1, \cdots, \alpha_k$. (Ces
opérations commutent puisque les arcs géodésiques sont disjoints.)

\subsection{Déformations expansives et complexe des arcs}

Soit $\rho:\pi_1S\to PSL(2,\R)$ la représentation d'holonomie de $S$,
dont on suppose qu'elle est convexe co-compacte. Pour tout $\gamma\in \pi_1S$,
on appelle $\lambda_\gamma(\rho)$ la longueur de translation hyperbolique de 
$\rho(\gamma)$.

\begin{defi}
Soit $j:\pi_1S\to PSL(2,\R)$ une autre représentation. On dit que $j$~est
{\em uniformément plus grande} que $\rho$ s'il existe $\epsilon>0$ tel que,
pour tout $\gamma\in \pi_1S\setminus \{ e\}$, 
$\lambda_\gamma(j)\geq (1+\epsilon)\lambda_\gamma(\rho)$.   
\end{defi}

\begin{defi}
On note $\Adm^+(\rho)$ l'espace des représentations de $j:\pi_1S\to PSL(2,\R)$ 
(considérées à conjugaison près) convexe co-compactes 
qui sont {\og uniformément plus grandes\fg} que $\rho$.
\end{defi}

Cet espace de représentations {\og plus grandes\fg} que $\rho$ admet une 
paramétrisation simple en termes de déformations par bandelettes. Cette
paramétrisation utilise le complexe des arcs, qui est l'un des complexes
simpliciaux qu'on peut associer à une surface compacte à bord (non vide).

\begin{defi}
Le {\em complexe des arcs}
de $S$ est le complexe simplicial $\bar X$ dont les $k$-faces sont déterminées
par les ensembles de $k+1$ arcs (topologiques) homotopiquement disjoints
de $S$. On note $X$ le sous-complexe constitué des faces pour 
lesquelles le complémentaire des arcs est une réunion disjointe de disques.
\end{defi}

Ainsi, si un point $t$ du complexe des courbes est contenu dans une face de 
dimension~$k$ correspondant à des arcs $\bar \alpha_1,\cdots,\bar \alpha_{k+1}$ 
homotopiquement disjoints, alors $t$ correspond à des poids 
$(t_1,\cdots, t_{k+1})$ positifs ou nuls, de somme égale à $1$, sur les $k+1$ arcs. 
Ces faces simpliciales sont recollées de la manière naturelle.

Il est aussi nécessaire de considérer le cône $CX$ sur $X$, c'est-à-dire qu'on
considère des poids dont la somme n'est pas nécessairement égale à $1$. 

On choisit pour chaque arc topologique $\bar \alpha$ sur $S$ un représentant 
géodésique $\alpha$, de manière que ces représentants s'intersectent de
manière minimale (c'est possible), ainsi qu'un point $p_\alpha\in \alpha$ et
une largeur $w_\alpha>0$. On peut maintenant énoncer le résultat suivant
\cite[Theorem 1.8]{danciger-gueritaud-kassel:margulis}.

\begin{theo}
L'application de $CX$ dans $\Adm^+(\rho)$ qui à une famille de poids
$t_{\alpha_i}$ sur des arcs homotopiquement disjoints $\bar \alpha_i$, $i=1,\cdots,k+1$ 
associe l'image de $\rho$ par la déformation en bandelettes
associées aux représentants géodésiques $\alpha_i$, aux points 
$p_{\alpha_i}$ et aux poids $t_{\alpha_i}w_{\alpha_i}$, est un homéomorphisme.
\end{theo}

Notons qu'on pourrait supprimer la dépendance aux $w_{\alpha_i}$, quitte à
prendre une paramétrisation légèrement différente. La dépendance par
rapport aux représentants géodésiques des arcs et aux points $p_\alpha$
est par contre bien réelle.

\subsection{Déformations par bandelettes infinitésimales}

On peut définir de même une version infinitésimale des déformations par
bandelettes. 

\begin{defi}
Soit $\alpha_1, \cdots, \alpha_k$ un ensemble d'arcs géodésiques disjoints,
munis de points $p_{\alpha_1}\in \alpha_1,\cdots,p_{\alpha_k}\in \alpha_k$ et de largeurs  
$w_{\alpha_1}, \cdots, w_{\alpha_k}$. Pour $t>0$, on note $\rho_t$ la déformation 
par bandelettes de $\rho$ associée aux poids $tw_{\alpha_1},\cdots, tw_{\alpha_k}$
sur les arcs géodésiques  $\alpha_1, \cdots, \alpha_k$ avec points marqués $p_1,
\cdots, p_k$. La déformation par bandelette infinitésimale associée à ces données
est $(d\rho_t/dt)_{|t=0}$, vu comme un vecteur tangent à l'espace des représentations
de $\pi_1S$ dans $PSL(2,\R)$.
\end{defi}

Comme pour les déformations {\og macroscopiques\fg}, ces déformations par bandelette
infinitésimale ne décroissent infinitésimalement la longueur géodésique d'aucune
courbe fermée. Notons $\adm^+(\rho)$ le cône des déformations infinitésimales
{\og expansives\fg}. Une fois choisis pour chaque arc (topologique) $\bar \alpha$ 
un représentant géodésique $\alpha$, un point $p_\alpha\in \alpha$, et
un poids $w_\alpha>0$, on a une application $f$ du complexe des arcs $\bar X$ dans
$\adm^+(\rho)$. On a alors la paramétrisation suivante 
\cite[Theorem 1.3]{danciger-gueritaud-kassel:margulis}.

\begin{theo} 
L'application $f$, composée avec la projectivisation $P:\adm^+(\rho)\to \adm^+(\rho)/\R_{>0}$,
définit un homéomorphisme entre $X$ et $\adm^+(\rho)/\R_{>0}$.
\end{theo}

\section{Variétés lorentziennes plates de dimension 3}

On va maintenant brièvement décrire quelques résultats concernant les espaces-temps
de Margulis, qui sont des variétés lorentziennes plates complètes dont le groupe 
fondamental est un groupe libre (non-abélien). La relation avec les déformations
expansives de métriques hyperboliques sur les surfaces apparaitra au $\S$2.4.

\subsection{Contexte : la conjecture d'Auslander et 
les espaces-temps de Margulis}

Un théorème classique de Bieberbach (voir par exemple \cite{auslander:theory})
décrit les groupes
fondamentaux des variétés compactes munies de métriques euclidiennes~: ce sont
des groupes discrets contenant un sous-groupe commutatif d'indice fini.

Il est naturel de se demander dans quelle mesure on peut étendre cet énoncé
d'une part à des structures géométriques moins rigides que des structures euclidiennes,
et d'autre part à des variétés qui ne sont pas nécessairement compactes. Un
{\em groupe cristallographique affine} est un groupe discret agissant 
proprement discontinuement
par transformations affines sur $\R^n$ avec quotient compact. On peut souhaiter
généraliser le théorème de Bieberbach et décrire ces groupes cristallographiques
affines, voire les variétés munies de structures affines complètes.

La conjecture d'Auslander, connue jusqu'en dimension 6 \cite{abels-margulis-soifer:6}, 
affirme que le groupe fondamental d'une variété affine
compacte complète est virtuellement poly\-cyclique. Pour les variétés seulement
complètes, la situation est plus riche. En 1977, Milnor \cite{milnor:fundamental-groups} avait demandé 
si un groupe libre non-abélien pouvait agir proprement discontinuement par transformation affines
sur $\R^n$, et suggéré de considérer un groupe libre agissant sur l'espace de
Minkowski de dimension $3$, $\R^{2,1}$. Un peu plus tard, en 1983, Fried et Goldman 
\cite{fried-goldman} montraient
que parmi les variétés affines complètes de dimension $3$, seules ces variétés
lorentziennes plates peuvent avoir un groupe fondamental qui n'est pas virtuellement
résoluble.

Enfin, toujours en 1983, Margulis \cite{margulis:free,margulis:complete} 
construisit des actions propre libres de
groupes libres non-abéliens sur $\R^{2,1}$, répondant ainsi à la question de Milnor.
On va voir ci-dessous comment caractériser les actions de ce type.

\subsection{Le groupe d'isométries de $\R^{2,1}$}

Le groupe des isométries de $\R^{2,1}$ est le produit semi-direct $O(2,1)\ltimes \R^{2,1}$,
où $O(2,1)$ agit en préservant l'origine, alors que $\R^{2,1}$ agit par translation.

Soit maintenant $\Gamma$ le groupe fondamental d'une surface, et soit
$\rho:\Gamma\to \isom(\R^{2,1})$. Alors la composée $j$ de $\rho$ avec
la projection $O(2,1)\ltimes \R^{2,1}$ sur le premier facteur, appelée partie linéaire de $\rho$, est un
morphisme de $\Gamma$ dans $O(2,1)$. On va considérer des morphismes 
dont la partie linéaire est à valeur dans la composante de l'identité
de $O(2,1)$, qu'on identifie avec $PSL(2,\R)$.

Par contre, si on note $\tau:\Gamma\to \R^{2,1}$
la composée de $\rho$ avec la projection sur le second facteur, $\tau$ est un 
1-cocycle pour $j$, c'est-à-dire qu'il satisfait la propriété d'équivariance
suivante~:
$$ \forall \gamma,\gamma'\in \Gamma, ~ \tau(\gamma.\gamma')=\tau(\gamma)+Ad(j(\gamma))\tau(\gamma') ~.$$

\subsection{La partie translation des cocycles comme une déformation}

Comme $\tau$ est un 1-cocycle pour $j$, on peut l'interpréter comme une déformation
infini\-tésimale de $j$, c'est-à-dire comme un vecteur tangent à l'espace des 
homéomorphismes de~$\Gamma$ dans $PSL(2,\R)$. 

Or Fried et Goldman \cite{fried-goldman}
ont montré que si $(j,\tau)$ est une action propre sur $\R^{2,1}$ d'un groupe non
virtuellement résoluble, alors $j(\gamma)$ doit être virtuellement un groupe de surface.
Mais $j(\Gamma)$ ne peut pas être co-compact.

Pour trouver des groupes agissant proprement sur $\R^{2,1}$ non virtuellement
résolubles, on est donc conduit à considérer une surface hyperbolique $S$
complète non compacte, munie d'une déformation infinitésimale $\tau$ de la métrique
(ou de manière équivalente de sa représentation d'holonomie). 

\subsection{La conjecture de Margulis}

Dans ce contexte, Margulis \cite{margulis:free} avait défini pour chaque élément 
$\gamma\in \Gamma$ un invariant $\alpha_ \tau(\gamma)$ (dépendant de $j$ et de $\tau$) et conjecturé que
la positivité (ou négativité) de cet invariant pour tous les éléments de $\Gamma$ assure 
la propreté de l'action sur $\R^{2,1}$. 

Goldman et Margulis \cite{goldman-margulis} ont ensuite donné une interprétation
simple de cet invariant~: c'est la variation infinitésimale, sous la 
déformation $\tau$, de la longueur hyper\-bolique de la géodésique fermée qui réalise 
$\gamma$. Le critère de Margulis est donc que, quitte à remplacer $\tau$ par $-\tau$, 
$\tau$ fait décro\^\i tre infinitésimalement la longueur de toutes les géodésiques
fermées sur $S$.

Goldman, Labourie et Margulis \cite{goldman-labourie-margulis} ont ensuite 
démontré une version légèrement précisée de cette conjecture~: l'action
$(j,\tau)$ de $\Gamma$ sur $\R^{2,1}$ est propre si et seulement si, quitte à
remplacer $\tau$ par $-\tau$, la 
déformation infinitésimale associée à $\tau$ décro\^\i t uniformément la longueur
de toutes les géodésiques fermées sur $S$, c'est-à-dire si et seulement si il
existe $\epsilon>0$ tel que, pour tout $\gamma\in \Gamma$, $\gamma\neq 1$, on a
$$ \frac d{dt}\left(
\frac{\lambda(e^{t\tau}j(\gamma))}{\lambda(j(\gamma))}
\right)_{|t=0}\leq -\epsilon~, $$
où $\lambda$ désigne la longueur de translation.

\subsection{Actions propres et champs de vecteurs contractants}

Danciger, Guéritaud et Kassel \cite[Theorem 1.1 (1)]{danciger-gueritaud-kassel:geometry}
donnent un autre critère, équivalent, pour la propreté de l'action $(j,\tau)$.

Soit $v$ un champ de vecteurs sur le plan hyperbolique $\HH^2$. On dira que 
$v$ est \mbox{{\em $(j,\tau)$-équivariant}} si~: 
$$ \forall \gamma\in \Gamma,\forall x\in \HH^2, 
v(j(\gamma)(x)) = (j(\gamma))_*(v(x)) + \tau(\gamma)(j(\gamma)(x))~. $$
Ici $j(\gamma)\in PSL(2,\R)$ est considéré comme une isométrie de $\HH^2$,
alors que $\tau(\gamma)\in sl(2,\R)$ est vu comme un champ de Killing sur $\HH^2$. 
Un champ de vecteurs $(j,\tau)$-équivariant induit donc $\tau$ comme variation de
la représentation $j$.

On dira par ailleurs qu'un champ de vecteurs est {\em $\epsilon$-contractant} si,
lorsqu'on note $h$ la métrique hyperbolique sur $\HH^2$, sa dérivée de Lie vérifie
$$ \cL_vh\leq -\epsilon h~, $$
c'est-à-dire que le flot de $v$ fait décro\^\i tre uniformément les longueurs des
vecteurs \mbox{tangents.}

Danciger, Guéritaud et Kassel \cite[Theorem 1.1]{danciger-gueritaud-kassel:geometry}
démontrent la version suivante de la conjecture de Margulis.

\begin{theo} \label{tm:proprete-minko}
Soient $\Gamma$ un groupe discret, et $(j,\tau):\Gamma\to PSL(2,\R)\ltimes sl(2,\R)$
un morphisme dont la partie linéaire $j$ est convexe co-compacte.
L'action de $(j,\tau)$ sur $\R^{2,1}$ est proprement discontinue si et seulement si,
quitte à remplacer $\tau$ par $-\tau$,  
il existe un champ de vecteurs $\epsilon$-contractant et $(j,\tau)$-équivariant
de $\HH^2$.
\end{theo}

\subsection{La topologie des quotients}

L'approche suivie par Danciger, Guéritaud et Kassel permet aussi de comprendre la
topologie des espaces-temps Minkowski qui peuvent être obtenus. En conséquence,
ils montrent un résultat de {\em sagesse} conjecturé par Drumm et Goldman
\cite{drumm-goldman}~: ces variétés sont homéomorphes à l'intérieur de variétés
compactes à bord. 

\begin{theo}[\cite{danciger-gueritaud-kassel:geometry}]
Soit $\Gamma$ un groupe discret sans torsion, et soit \mbox{$j\in \Hom(\Gamma,
PSL(2,\R))$} une représentation convexe co-compacte dont le quotient est
une surface $S$. Soit $\tau$ un cocycle de déformation de $j$ tel que
$(j,\tau)$ agit proprement sur~$\R^{2,1}$. Alors la variété quotient 
est un fibré en droite sur $S$, avec pour fibres des géodésiques de
type temps. En conséquence, cette variété quotient est homéomorphe
à l'intérieur d'un corps à anse.
\end{theo}

On peut noter qu'une autre démonstration de la conjecture de sagesse
a été annoncée par Choi et Goldman, voir \cite{choi-goldman:tameness}.

On verra dans la section \ref{sc:geodesiques} un lien entre 
applications contractantes et feuilletages par des géodésiques
de type temps, qui permet de mieux comprendre cet énoncé.

\subsection{Domaines fondamentaux et plans croches}

Les plans croches ont été introduits par Drumm \cite{drumm} pour 
construire des actions propres de groupes libres sur $\R^{2,1}$ en
construisant directement des domaines fondamentaux pour des 
actions dont la partie linéaire est la représentation d'holonomie
d'une surface hyperbolique complète. 

Un plan croche dans $\R^{2,1}$ est centré en un point $c$, et composé de trois parties~:
\begin{itemize}
\item une partie centrale, qui est l'intersection du cone de 
lumière de $c$ avec un plan de type temps passant par $c$,
\item deux demi-plans isotropes (de type lumière) dont les
bords sont recollés avec les deux droites qui constituent le
bord de la partie centrale. 
\end{itemize}

\begin{figure}[h]
\begin{center}
\includegraphics[width = 6.0cm]{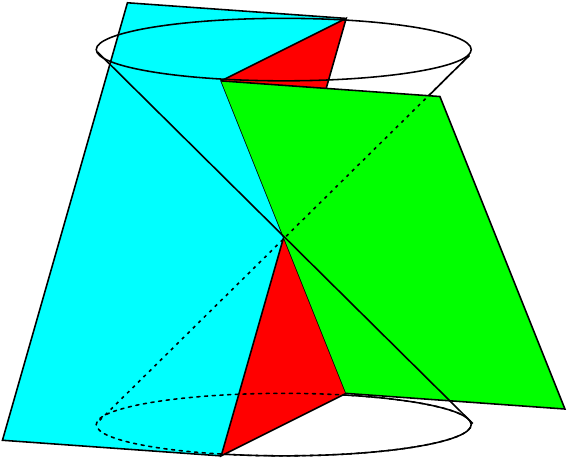}
\caption{Un plan croche}
\end{center}
\end{figure}

Il existe deux types de plan croche, les {\og gauches\fg} et les {\og droits\fg},
mais pour construire des domaines fondamentaux, on n'utilise que l'un
des deux types, par exemple seulement des plans croches gauches.

Un plan croche gauche de $\R^{2,1}$ est uniquement déterminé
par sa partie centrale, qui correspond, à translation près, 
à une géodésique dans $\HH^2$. Si on se
donne une famille de géodésiques disjointes dans $\HH^2$ (par exemple
invariante sous une action convexe co-compacte de $\Gamma$ sur $\HH^2$), 
on peut donc associer un plan croche défini à translation près à chaque géodésique, 
mais la position du centre doit être déterminée par la composante
translation $\tau$ de l'action de $\Gamma$ sur $\R^{2,1}$. L'une des raisons
pour lesquelles les plans croches sont intéressants et utiles est qu'on peut
déterminer précisément les conditions sur les positions des centres sous
lesquelles deux plans croches gauches sont disjoints, ce qui a permis
à Drumm d'obtenir les résultats mentionnés plus haut.

Ces constructions ont conduit Drumm et Goldman à conjecturer que 
les espaces-temps de Margulis admettent toujours un domaine
fondamental délimité par des plans croches, conjecture démontrée
par Charette, Drumm et Goldman pour les groupes libres de rang~2
\cite{charette-drumm-goldman1,charette-drumm-goldman2}. 

Danciger, Guéritaud et Kassel démontrent cette conjecture dans le cas général
\cite[Theorem 1.7]{danciger-gueritaud-kassel:margulis}.

\begin{theo} \label{tm:croches-minko}
Soit $\Gamma$ un sous-groupe discret de $\isom(\R^{2,1})$ agissant
proprement discontinuement et librement sur $\R^{2,1}$ avec pour
partie linéaire une représentation convexe co-compacte. Alors 
l'action de $\Gamma$ sur $\R^{2,1}$ admet un domaine fondamental
bordé par un nombre fini de plans croches.
\end{theo}

On peut déjà noter à ce stade que leur preuve est une conséquence
assez directe de la paramétrisation des déformations expansives par
les déformations en bandelettes infinitésimales --- on développera
ce point dans le $\S$\ref{sc:defos-croches}.

\section{Espace-temps anti-de Sitter}

On passe maintenant de l'espace de Minkowski à l'espace anti-de Sitter,
dont on va rappeler rapidement la définition et quelques propriétés essentielles.

\subsection{L'espace anti-de Sitter et $PSL(2, \R)$}
\label{ssc:ads}

On note $\R^{2,2}$ l'espace vectoriel $\R^4$ muni de la forme 
bilinéaire symétrique de signature $(2,2)$~:
$$ \langle x,x'\rangle_{2,2} = x_1x'_1+x_2x'_2-x_3w'_3-x_4x'_4~. $$ 
L'espace anti-de Sitter de dimension 3 peut se définir comme la quadrique~:
$$ AdS_3 = \{ x\in \R^{2,2}~|~ \langle x,x\rangle =-1\}~, $$
munie de la métrique induite. Il apparaît ainsi comme un analogue
lorentzien de l'espace hyperbolique de dimension 3~: sa courbure
sectionnelle est constante égale à $-1$, et ses plans totalement géodésiques
de type espace sont isométriques au plan hyperbolique. On note que
les plans totalement géodésiques de $AdS_3$ sont les intersections avec
$AdS_3$, vu comme une quadrique de $\R^{2,2}$, des hyperplans contenant
$0$. 

On peut aussi identifier $AdS_3$ avec $PSL(2,\R)$ muni de sa métrique
de Killing bi-invariante. On a une action naturelle de 
$PSL(2,\R)\times PSL(2,\R)$ sur $AdS_3$, 
les deux facteurs agissant par multiplication respectivement à gauche et à 
droite. Vu ainsi, $AdS_3$ rappelle la sphère $S^3$, qu'on peut 
identifier avec $SU(2)$ muni de sa métrique de Killing bi-invariante
et dont le groupe d'isométrie s'identifie (à indice fini près) à
$O(3)\times O(3)$.\footnote{Certaines des identifications ci-dessous devront
aussi s'entendre à indice 2 près.}

Enfin on peut voir $AdS_3$ comme l'espace total du fibré unitaire tangent au
plan hyperbolique, muni de sa métrique lorentzienne naturelle. On voit ainsi qu'il
existe beaucoup de variétés AdS fermées de dimension 3~: les fibrés tangents 
unitaires des surfaces hyperboliques fermées. On peut vérifier que les
représentations d'holonomie de ces variétés s'écrivent, dans la décomposition 
$\isom_0(AdS_3)=PSL(2,\R)\times PSL(2,\R)$, sous la forme $(\rho, 1)$, où 
$\rho$ est la représentation d'holonomie de la surface hyperbolique 
sous-jacente et $1$ est la représentation triviale.

On dispose d'un modèle projectif de $AdS_3$, ou plus précisément d'un
hémisphère $AdS_{3,+}$ de $AdS_3$, c'est-à-dire d'un domaine fondamental pour l'action 
de $\Z/2\Z$ qui envoie un point $x\in AdS_3$ sur $-x$, dans la définition
ci-dessus où $AdS_3$ est identifié à une quadrique de $\R^{2,2}$. Ce modèle
est projectif au sens où il envoie $AdS_{3,+}$ sur l'intérieur d'un 
hyperboloïde à une nappe dans $\R^3$, et envoie les géodésiques de 
$AdS_{3,+}$ sur les segments de droite. 

Il est analogue au modèle de
Klein de l'espace hyperbolique, et obtenu de la même manière~: 
on note $H_0$ le plan d'équation $x_3=1$, tangent à $AdS_3$ en $(0,0,1,0)$, 
et on projette tous les points de $AdS_3$ pour lesquels $x_3>0$ sur $H_0$
dans la direction de l'origine. Ainsi $AdS_{3,+}$ est envoyé sur l'intérieur
d'un hyperboloïde à une nappe. On peut enrichir ce modèle en le voyant non
pas dans $\R^3$ mais dans $\R P^3$ --- on obtient ainsi un modèle projectif
de $AdS_3/(\Z/2\Z)$ --- ou dans son revêtement à deux feuillets $S^3$ --- 
on a alors un modèle projectif de tout $AdS_3$.

On utilisera aussi plus bas une notion de dualité entre plans orientés de
type espace et points dans $AdS_3$. Soit $H\subset AdS_3$ un plan de type
espace, alors $H=AdS_3\subset H_0$, où $H_0\subset \R^{2,2}$ est un 
hyperplan dont la métrique induite est de signature $(2,1)$. Si $H$
est orienté, on en déduit une orientation de $H_0$. Le point $H^*$ dual de 
$H$ est défini comme la normale unitaire orientée à $H_0$, qui est 
bien dans $AdS_3$. On vérifie facilement que $H^*$ est contenu dans 
toutes les géodésiques orthogonales à $H$, et qu'il se trouve à distance 
$\pi/2$ de $H$ le long de chacune de ces géodésiques.

\subsection{Actions propres sur $AdS_3$}

On comprend maintenant bien les variétés AdS fermées de dimension 3, ou du moins
leurs représentations d'holonomie. 

Kulkarni et Raymond  \cite{kulkarni-raymond}
avaient montré que si $(j,\rho)$ agit proprement discontinuement
sur $AdS_3$ sans torsion, alors (à échange des deux facteurs près)
$j$ est injective et discrète. Plus récemment, ce résultat a été 
précisé par Kassel \cite{kassel:2012}, et on a le~:

\begin{theo}[Kassel] \label{tm:proprete-ads}
Si $j$  est convexe co-compacte, $(j,\rho)$ agit proprement
si et seulement si l'une des deux conditions équivalentes est vérifiée~:
\begin{itemize}
\item il existe une application $f:\HH^2\to \HH^2$ $(j,\rho)$-équivariante
qui est $k$-Lipschitz pour un $k<1$,
\item les longueurs de translation des éléments de $\Gamma$ sont 
uniformément plus courtes pour $\rho$ que pour $j$, au sens où
$$ \sup_{\gamma\in \Gamma, \lambda(j(\gamma))>0}\frac{\lambda(\rho(\gamma))}{\lambda(j(\gamma))}<1~. $$
\end{itemize}
\end{theo}

On peut voir cet énoncé comme l'analogue \og macroscopique\fg du théorème 
\ref{tm:proprete-minko} ci-dessus.

\subsection{Topologie et sagesse des variétés anti-de Sitter}

Comme pour le cas lorentzien plat décrit ci-dessus, Danciger, Guéritaud et Kassel
décrivent la topologie des variétés anti-de Sitter compl\`etes basées sur une 
représentation $j$ convexe co-compacte, voir 
\cite[Theorem 1.2 (1)]{danciger-gueritaud-kassel:margulis}.

\begin{theo}
Soit $\Gamma$ un groupe discret sans torsion, et soit $j$ une représentation 
convexe co-compacte de $\Gamma$ dans $PSL(2,\R)$ dont on note $S$ la surface quotient.
Soit $\rho$~une autre représentation de $\Gamma$ dans $PSL(2,\R)$ telle que $(j,\rho)$
agit proprement sur~$AdS_3$. Alors le quotient est un fibré en cercle dont les fibres
sont des géodésiques de type temps.
\end{theo}

On en déduit que ces variétés anti-de Sitter sont des fibrés de Seifert sur un
orbifold hyperbolique.

De plus, Danciger, Guéritaud et Kassel \cite[Theorem 1.9]{danciger-gueritaud-kassel:margulis}
montrent que les variétés AdS complètes
obtenues à partir de deux représentations convexe co-compactes admettent
un {\og bon\fg} domaine fondamental, comme on l'a vu dans le théorème \ref{tm:croches-minko} 
pour le cas Minkowski.

La notion de plan croche dans $AdS_3$ est similaire à celle vue plus haut dans
l'espace de Minkowski. Un plan croche est donc constitué de trois parties~:
\begin{itemize}
\item la partie centrale, qui est l'intersection du cone de lumière
d'un point $c$ avec un plan de type temps contenant $c$,
\item deux demi-plans de type lumière dont les bords co\"\i ncident avec
les deux composantes connexes du bord de la partie centrale privée de $c$.
\end{itemize}

\begin{figure}[h]
\begin{center}
\includegraphics[width = 8.0cm]{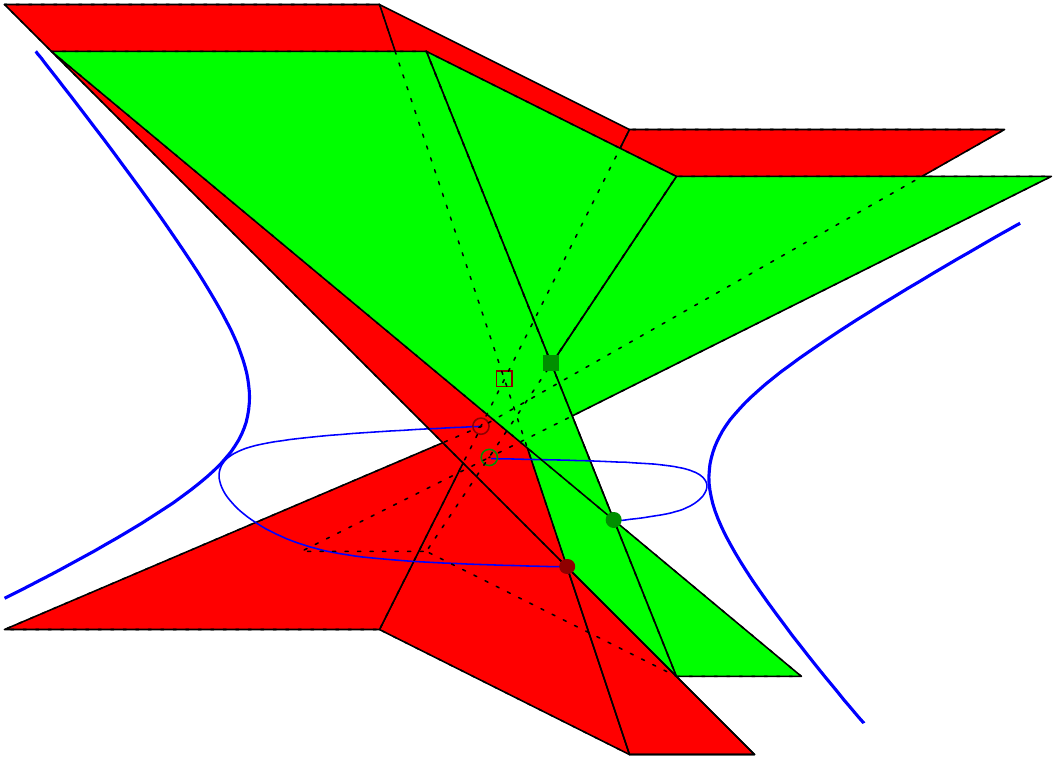}
\caption{Deux plans croches disjoints dans $AdS_3$}
\end{center}
\end{figure}

On a le résultat suivant \cite[Theorem 1.9]{danciger-gueritaud-kassel:margulis}.

\begin{theo}[Danciger, Guéritaud et Kassel] 
\label{tm:croches-ads}
Soient $\rho$ et $j$ les représentations d'holonomie de deux structures
convexes co-compactes sur $S$, telles que $(\rho,j)$ agit proprement
discontinuement sur $AdS_3$. Alors cette action admet un domaine fondamental
bordé par des plans croches.   
\end{theo}

La preuve est encore une conséquence de l'existence de déformations
par bandelettes (non infinitésimales) des surfaces hyperboliques.

Goldman \cite{goldman:crooked} a remarqué que les plans croches tant
Minkowski que anti-de Sitter peuvent être vus comme des restrictions des
plans croches dans l'univers d'Einstein, développés par Frances \cite{frances:margulis}.

\section{De anti-de Sitter à Minkowski~: la géométrie transitionnelle}
\label{sc:transition}

On note une analogie assez forte entre la description des actions propres sur
l'espace de Minkowski, en termes de décroissance de la longueur, et sur $AdS_3$,
en termes d'application contractante équivariante. Cette analogie s'explique
par l'existence d'une {\em géométrie transitionnelle}, qui offre une certaine
continuité entre les géométries hyperbolique, Minkowski et anti-de Sitter de 
dimension trois. 

\subsection{Géométrie demi-lune et géométrie Minkowski}
\label{ssc:demi-lune}

La géométrie transitionnelle développée par Danciger fait apparaître un objet
d'apparence exotique mais assez naturel, la géométrie {\og demi-lune\fg}.\footnote{Le terme 
{\og demi-lune\fg} est utilisé par certains adeptes francophones de
la planche à roulettes pour traduire le terme anglais {\og half-pipe\fg}. On gardera
ici la notation $\HP^3$, qui vient du {\og half-pipe\fg} anglais, pour désigner cet
espace en dimension 3.}
L'espace demi-lune de dimension $3$, $\HP^3$, est simplement $\HH^2\times \R$, muni
de la métrique dégénérée $0\times dt^2 +h$, où $h$ est la métrique du plan 
hyperbolique. 

Le groupe des isométries de $\HP^3$ est de dimension infinie, puisque la
métrique ne {\og voit\fg} pas les translations verticales. Il existe pourtant 
un sous-groupe du groupe des isométries, qui joue le rôle de groupe
d'isométries {\og restreint\fg}. C'est le produit $PSL(2,\R)\ltimes sl(2,\R)$, 
où le facteur $PSL(2,\R)$ agit par {\og translations horizontales\fg} agissant
sur $\HH^2$, et les éléments de $sl(2,\R)$ agissent par translations 
verticales proportionelles au rotationnel du champ de Killing de $\HH^2$ 
qui leur est associé. On a aussi une notion naturelle de plans totalement
géodésiques, ce sont les images par les éléments du groupe d'isométrie
restreint du plan horizontal $\HH^2\times \{ 0\}$. C'est cet espace $\HP^3$
qui apparaît le plus directement comme quand on considère une suite 
d'actions sur $AdS_3$ qui dégénèrent d'une manière particulière.

L'espace demi-lune $\HP^3$ est dual de l'espace de Minkowski au sens 
suivant. \`A~chaque plan de type espace $P$ dans Minkowski, on peut déterminer
un point dans~$\HP^3$, déterminé par la normale unitaire à $P$ orientée vers le
futur (identifiée à un point de $\HH^2$) et la distance orientée de l'origine
à $P$. L'angle entre deux plans correspond alors à la distance dans $\HP^3$ 
entre les points duaux. Réciproquement, à chaque plan totalement géodésique
de $\HP^3$, on associe un point de $\R^{2,1}$, et ces deux relations de dualité
réciproques sont compatibles. (On a une dualité analogue entre l'espace euclidien
$\R^3$ et l'espace $S^2\times \R$ muni de la métrique dégénérée $h_{S^2}+0dt^2$,
où $h_{S^2}$ est la métrique canonique de la sphère.)

Cette dualité permet de considérer une action d'un groupe $\Gamma$ sur $\HP^3$
comme une action sur $\R^{2,1}$, à travers l'action de $\Gamma$ sur les plans
totalement géodésiques de $\HP^3$. 

\subsection{Trois points de vue sur la géométrie transitionnelle}

La notion de géométrie transitionnelle est développée dans la thèse de Jeff
Danciger \cite{danciger:transition,danciger:ideal}. 
Elle permet de décrire précisément une forme de 
continuité entre variétés hyperboliques, Minkowski et anti-de Sitter.
Trois points de vue complémentaires au moins sont possibles~: en termes
de représentation, de géométrie projective, ou de  métriques (riemanniennes
ou lorentziennes suivant le cas).

\subsubsection{Transitions entre représentations et nombres de Lorentz}

Le groupe des isométries préservant l'orientation du plan hyperbolique,
$SO(2,1)$, peut s'identifier avec $PSL(2,\R)$. 

Considérons maintenant au lieu de $\C=\R+i\R$ l'anneau $\R+\tau\R$, avec 
$\tau^2=1$. C'est un
anneau bien connu qui porte plusieurs noms distincts, on parle en particulier
de nombres de Lorentz ou de nombres hyperboliques. Posons $\omega_\pm
=\frac{1\pm \tau}2$. On vérifie immédiatement que $\omega_\pm^2=\omega_\pm$, 
alors que $\omega_+\omega_-=0$. Donc si $A_+,A_-\in M_2(\R)$, on a 
$$ \det(\omega_+A_++\omega_-A_-)=\omega_+\det(A_+)+\omega_-\det(A_-)~. $$
Ainsi, si $A=\omega_+A_++\omega_-A_-$, alors $\det(A)=1$ si et seulement
si \mbox{$\det(A_+)=\det(A_-)=1$.} En poursuivant ce raisonnement, on voit qu'on peut 
identifier $PSL(2,\R+\tau\R)$ à $PSL(2,\R)\times PSL(2,\R)$. 

De plus, 
si $A=\omega_+A_++\omega_-A_-$ et $B=\omega_+B_++\omega_-B_-$ sont deux matrices
de \mbox{$SL(2,\R+\tau\R)$,} avec $A_+,A_-,B_+,B_-\in SL(2,\R)$, alors 
$$ AB=(\omega_+A_++\omega_-A_-)(\omega_+B_++\omega_-B_-) =
\omega_+ A_+B_++\omega_-A_-B_-~. $$
La structure produit de  $PSL(2,\R+\tau\R)$ correspond donc à celle de
$PSL(2,\R)\times PSL(2,\R)$. 

On peut donc identifier $\isom_0(AdS_3)$ avec $PSL(\R+\tau\R)$. 
On note au passage que l'injection canonique
de $PSL(2,\R)$ dans $PSL(2,\R+\tau\R)$ a pour image la diagonale dans 
$PSL(2,\R)\times PSL(2,\R)$.

Remplaçons maintenant $\R+\tau \R$ par $\R+\epsilon \R$, avec $\epsilon^2=0$. Soit 
$A_0, A_1\in M_2(\R)$, alors 
$$ \det(A_0(I+\epsilon A_1)) = \det(A_0)(1+\epsilon \mathrm{tr}(A_1))~, $$
et donc $\det(A_0(I+\epsilon A_1))=1$ si et seulement si $\det(A_0)=1$
et $\mathrm{tr}(A_1)=0$. On peut donc identifier $SL(2,\R+\epsilon\R)$ avec 
$SL(2,\R) \times sl(2,\R)$. 

De plus, si $A_0(I+\epsilon A_1), B_0(I+\epsilon B_1)\in SL(\R+\epsilon \R)$, 
alors
$$ A_0(I+\epsilon A_1) B_0(I+\epsilon B_1) = 
A_0B_0 + \epsilon(A_0A_1B_0 + A_0B_0B_1) =
A_0B_0(I+\epsilon(B_0^{-1}A_1B_0+B_1))~, $$
et on retrouve la multiplication du produit semi-direct. 
En poursuivant, on identifie ainsi $\isom_0(\R^{2,1})$ avec
$PSL(2,\R)\ltimes sl(2,\R)$. 

Considérons maintenant une famille régulière à un paramètre $(j(t), \rho(t))$
de morphismes de $\Gamma$ dans $PSL(2,\R)\times PSL(2,\R)$, où 
$\Gamma$ est le groupe fondamental d'une surface $S$ (à bord),
et supposons que $\rho(0)=j(0)$. Posons $u(t)=\omega_+j(t)+\omega_-(t)\rho$,
on a alors 
$$ u(t)=(j(t)+\rho(t))/2 + \tau(j(t)-\rho(t))/2~. $$
Si on pose $v=(j(t)-\rho(t))'(0)/2$, alors $v$ est un 1-cocycle
pour $u(0)$, et $u(0)+\epsilon v$ détermine
une représentation de $\Gamma$ dans $PSL(2,\R+\epsilon\R)$. Inversement,
si on se donne une représentation $u$ de $\Gamma$ dans $PSL(2,\R)$ et 
un 1-cocycle $v$ pour $u$, on obtient pour $t$ petit une paire de 
représentations $u\exp(-tv)$ et $u\exp(tv)$ à valeurs dans $PSL(2,\R)$, 
et donc une représentation de $\Gamma$ dans $\isom_0(AdS_3)$.

Un point essentiel est que la condition pour que $(j(t),\rho(t))$ agisse
proprement sur $AdS_3$ est que (quitte à échanger les facteurs) les longueurs
de translation de tous les éléments de $\Gamma$ soient uniformément plus
petites pour $\rho$ que pour $j$. Pour les représentations dans $\isom_0(\R^{2,1})$,
la condition correspondante est que le cocycle de déformation raccourcisse
uniformément les longueurs de translation de tous les éléments de $\Gamma$. 
Les deux conditions correspondent parfaitement, et la {\og dérivée en $0$\fg} d'une
famille à un paramètre d'actions propres sur $AdS_3$, telle que les deux
représentations coïncident à la limite $t\to 0$, fournit une action propre
sur $\R^{2,1}$, et réciproquement.

On voit dans cette description que, dans la représentation à valeur dans 
$PSL(2,\R)\ltimes sl(2,\R)$ obtenue à la limite, la composante dans 
$PSL(2,\R)$ correspond à la limite en $0$ des deux représentations 
$j(t)$ et $\rho(t)$, alors que le cocycle à valeurs dans $sl(2,\R)$ décrit
les différences entre les dérivées en $0$ de $j$ et $\rho$. La même 
interprétation est apparente dans les autres descriptions données ci-dessous,
on n'y reviendra pas.

\subsubsection{Transitions projectives}

On peut aussi donner un point de vue projectif sur les transitions
géométriques. Soit encore $(j(t),\rho(t))_{t\in ]0,1]}$ une famille à
un paramètre d'actions de $\Gamma$ sur $AdS_3$. Pour $t\in ]0,1]$,
on peut voir $(j(t),\rho(t))$ comme une action de $\Gamma$ sur
$\R P^3$ qui préserve une quadrique de signature $(1,1)$, qui n'est
autre que le bord à l'infini de $AdS_3$ dans le modèle projectif décrit
au $\S$\ref{ssc:ads}. 
 
Supposons que pour $t\to 0$, $(j(t),\rho(t))$ stabilise un plan
de type espace $H\subset AdS_3$. Alors l'action limite 
$(j(0),\rho(0))$ sur $\R P^3$ stabilise $\partial AdS_3$ ainsi
que $\partial P $ (le bord à l'infini de $P$ vu comme sous-ensemble
de $\partial AdS_3$) et il stabilise donc l'enveloppe des plans 
tangents à $\partial AdS_3$ le long de $\partial P$, enveloppe 
qui borde un {\og cylindre\fg}~$C$. (Ce {\og cylindre\fg} $C$
n'est autre, en termes de géométrie AdS, que le domaine de dépendance
de $P$.)

Or $C$ est naturellement un modèle projectif de $\HP^3$, on peut le
voir en considérant la métrique de Hilbert de $C$, qui est bien 
une métrique dégénérée dont la restriction aux intersections de $C$
avec presque tous les plans est une métrique hyperbolique.

On obtient ainsi pour $t=0$ une action de $\Gamma$ sur $\HP^3$, 
qui peut être vue comme une action sur $\R^{2,1}$ grâce à la dualité
mentionnée au $\S$\ref{ssc:demi-lune}.

On note que cette description s'applique de la même manière pour une
famille d'actions sur l'espace hyperbolique $\HH^3$, lorsque la limite
laisse stable un plan totalement géodésique.

\subsubsection{Transitions entre métriques}

On se place ici encore dans le cas où, pour $t\to 0$, $(j(t),\rho(t))$ stabilise un plan
de type espace $H\subset AdS_3$, et converge vers une action sur~$H$. On peut
considérer une suite d'applications du domaine de dépendance de $H$ à valeurs dans
$H\times ]-f(t),f(t)[$, avec $f(t)\to \infty$ quand $t\to 0$, qui envoie sur 
les lignes verticales $\{x\}\times ]-f(t),f(t)[$ les segments géodésiques de
longueur $\pi$ orthogonaux à $H$ dans $AdS_3$. 
Si on choisit convenablement la normalisation dans le facteur vertical, on peut
obtenir à la limite une action de $\Gamma$ sur $\HP^3$.

\subsubsection{Autres points de vue}

D'autres points de vue encore sont possibles sur la géométrie transitionnelle.
Danciger \cite{danciger:ideal} utilise des triangulations
idéales et les paramétrisations des espaces de représentations qui leurs sont
associées. Un autre point de vue encore, utilisé par exemple dans 
\cite[Appendix 2]{earthquakes}, consiste à {\og zoomer\fg} sur le point dual dans
$AdS_3$ (resp. dans l'espace de Sitter $dS_3$) du plan totalement géodésique laissé 
invariant par la limite d'une famille
d'actions sur $AdS_3$ (resp. $\HH^3$), pour voir apparaître une action sur $\R^{2,1}$ 
et par dualité sur $\HP^3$.

\subsection{Régénération d'espace-temps de Margulis en actions sur $AdS_3$}

En utilisant ces notions de géométrie transitionnelles, Danciger, Guéritaud et Kassel
\cite[Theorem 1.4]{danciger-gueritaud-kassel:geometry} décrivent comment la
représentation d'holonomie d'un espace-temps de Margulis peut se {\og régénérer\fg}
en une action propre sur $AdS_3$.

\begin{theo}
Soit $M=(j,u)(\Gamma)\backslash \R^{2,1}$ un espace-temps de Margulis tel que
\mbox{$S=j(\Gamma)\backslash \HH^2$} est une surface hyperbolique convexe co-compacte. 
Soient $t\to j_t$ et $t\to \rho_t$ des applications régulières telles que 
$j_0=\rho_0=j$ et que $\frac d{dt}|_{t=0}\rho_tj_t^{-1}=u$. Alors~:
\begin{itemize}
\item pour tout $t>0$ assez petit, $(j_t,\rho_t)(\Gamma)$ agit 
proprement discontinuement sur $AdS_3$,
\item  il existe une famille régulière de difféomorphismes $(j_t,\rho_t)$-invariants
de \mbox{$H^2\times S^1$} dans $AdS_3$, définis pour $t>0$ assez petit, déterminant des 
structures AdS complètes $\cA_t$ sur la variété fixée $S\times S^1$,
\item les structures projectives réelles $\cP_t$ sous-jacentes à $\cA_t$ convergent 
quand $t\to 0$ vers une structure projective réelle sur $S\times S^1$, et l'espace-temps
de Margulis $M$ est la restriction de $\cP_0$ à $S\times ]-\pi,\pi[$.
\end{itemize}
\end{theo}

Danciger, Guéritaud et Kassel en déduisent un autre énoncé plus {\og métrique\fg} qui met en évidence
le rôle des géodésiques de type temps dans la convergence 
\cite[Corollary 1.5]{danciger-gueritaud-kassel:geometry} 

\begin{coro}
Sous les hypothèses du théorème précédent, il existe une famille de métriques
AdS complètes $(g_t)_{t\in [0,\epsilon]}$ sur $S\times S^1$ telles que les
restrictions des $g_t$ à $S\times ]-\pi,\pi[$ convergeant uniformément sur les
compacts vers une métrique lorentzienne plate $g$ telle que $(S\times ]-\pi,\pi[,g)$
est isométrique à $M$.
\end{coro}

\section{Géodésiques de type temps}
\label{sc:geodesiques}

Dans cette section, on introduit quelques propriétés des géodésiques de 
type temps dans les variétés anti-de Sitter, qui permettent de mieux
comprendre les critères de propreté des théorèmes \ref{tm:proprete-ads} 
et \ref{tm:proprete-minko}, ainsi que l'existence
d'un feuilletage par des géodésiques de type temps pour les quotients
de $AdS_3$ par les groupes agissant proprement.

\subsection{Géodésiques de type temps dans $AdS_3$}
\label{ssc:temps-ads}

On a vu qu'on peut identifier $AdS_3$ avec $PSL(2,\R)$ muni de sa métrique
de Killing, lorentzienne et bi-invariante. Dans cette identification, on
vérifie sans difficulté que l'ensemble des rotations de $\R^2$ centrées en $0$
forme une géodésique de type temps.
Si on identifie $PSL(2,\R)$ avec les isométries du plan hyperbolique, 
on en déduit que chaque géodésique de type temps correspond à l'ensemble
des isométries qui envoient $x$ sur $y$, où $x$ et $y$ sont deux points
donnés de $\HH^2$. On peut donc identifier l'espace $\cG$ des géodésiques de 
type temps de $AdS_3$ avec $\HH^2\times \HH^2$.

Un autre point de vue, plus local, est possible (et utilisé dans \cite{colII}).
On peut définir deux connexions $D^l, D^r$ sur le fibré des vecteurs unitaires de 
type temps de $AdS_3$ de la manière suivante. Soit $n$ un champ de vecteurs unitaire
de type temps, et soit $x$ un vecteur quelconque~; on pose~:
$$ D^l_xn = \nabla_xn + x\times n~, ~~D^r_xn = \nabla_xn - x\times n~, $$
où $\nabla$ est la connexion de Levi-Civita de $AdS_3$ et $\times$ désigne le
produit vectoriel lorentzien. On vérifie alors par un calcul direct (voir \cite{colII})
que $D^l$ et $D^r$ sont des connexions plates et nulles le long du flot 
géodésique. On obtient ainsi deux projections de l'espace des géodésiques de type temps
de AdS sur l'espace des vecteurs unitaires de type temps en un point, lui-même 
identifié au plan hyperbolique.

Un troisième point de vue peut être utile, et sera utilisé dans le 
$\S$\ref{ssc:defos-croches-ads}. Il fournit une application explicite qui à
une géodésique $\delta\subset AdS_3$ de type temps associe un couple de
points dans $\HH^2$, identifié à un plan totalement géodésique fixé $H_0\subset AdS_3$. 
La construction repose sur le fait que $\partial_\infty AdS_3$, vu comme une
quadrique dans $\R^3$ par le modèle projectif décrit au $\S$\ref{ssc:ads}, 
est feuilleté par deux familles de droites. On parlera des feuilletages gauche
et droit. 

(Un autre point de vue encore est développé dans \cite{minsurf}, en termes de 
deux métriques hyperboliques définies sur des surfaces de type espace {\og pas
trop courbées\fg} dans une variété anti-de Sitter.)

\subsection{Feuilletages géodésiques et difféomorphismes de $\HH^2$}

Considérons une famille à un paramètre de géodésiques de type temps, 
$(c_t)_{t\in [0,1]}$. Soient $c_l(t)$ et $c_r(t)$ les projections sur les
facteurs dans l'identification de l'espace des géodésiques $\cG$ avec 
$\HH^2\times\HH^2$. Un calcul immédiat utilisant les connexions $D^l$ et $D^r$ 
introduites ci-dessus permet de montrer le~:

\begin{lemm} \label{lm:intersection}
Il existe $x\in c(0)$ dont la distance
à $c(t)$ a une dérivée nulle par rapport à $t$ en $t=0$ 
si et seulement si $\| c'_l(0)\|=\| c'_r(0)\|$.
\end{lemm}

Heuristiquement, cette propriété correspond au fait que $c(t)$ a une intersection avec 
$c(0)$ au premier ordre en $t=0$.

Cet énoncé a est la version infinitésimale du fait élémentaire suivant~: deux
géodésiques de type temps $c_0, c_1$ de $AdS_3$ se rencontrent si et seulement si les couples
de points de $\HH^2$ qui leurs sont associés, $(c_{0,l}, c_{0,r})$ et $(c_{1,l}, c_{1,r})$,
sont tels que $d(c_{0,l}, c_{0,r})=d(c_{1,l}, c_{1,r})$. C'est en effet à cette condition qu'il
existe une isométrie de $\HH^2$ qui envoie $c_{0,l}$ sur $c_{0,r}$ et $c_{1,l}$ sur $c_{1,r}$,
isométrie qui correspond au point d'intersection entre les deux géodésiques.

On considère maintenant un difféomorphisme $\phi:\HH^2\to \HH^2$. On peut voir son
graphe dans $\HH^2\times \HH^2$ comme une famille à deux paramètres de géodésiques
de type temps dans $AdS_3$, $\cF\subset \cG$. Une conséquence directe du lemme 
\ref{lm:intersection} est que $\cF$ est localement un feuilletage 
de $AdS_3$ dès que $\phi$ est strictement contractante
(ou dilatante), c'est-à-dire si la norme de sa différentielle est partout inférieure
à $k<1$. En poursuivant ce raisonnement, on peut montrer que, sous cette hypothèse,
$\cF$ est en fait un feuilletage de $AdS_3$.

Soit maintenant $(j,\rho)$ une action de $\Gamma$ sur $\HH^2\times \HH^2$ 
telle que $j$ est convexe co-compacte. Supposons que $\phi$ est équivariante
pour $(j,\rho)$, alors $\cF$ est invariant sous l'action de $\Gamma$, si 
bien que $\Gamma$ agit sur l'espace des feuilles de $\cF$.
Comme $j$ agit proprement, cette action est propre. On peut alors en 
déduire assez directement que l'action de $\Gamma$ sur $AdS_3$ est elle
aussi propre.

\subsection{Géodésiques de type temps dans $\R^{2,1}$}
\label{ssc:temps-minko}

Comme dans $AdS_3$, on peut identifier l'espace $\cG_0$ des géodésiques de
type temps dans $\R^{2,1}$ avec l'espace total $T\HH^2$ du fibré tangent au
plan hyperbolique. Considérons pour cela un point $(n,v)\in T\HH^2$, c'est-à-dire
que $n\in \HH^2$ et que $v\in T_n\HH^2$. On peut voir $n$ comme un vecteur 
unitaire de type temps orienté vers le futur dans $\R^{2,1}$, et $v$ comme un
vecteur de type espace orthogonal à $n$. On associe alors à $(n,v)$ l'unique
droite de type temps parallèle à $n$ et passant par $n\times v$, où
$\times$ désigne le produit vectoriel usuel de l'espace de Minkowski.

On a dans ce contexte un analogue direct du lemme \ref{lm:intersection}.
Considérons une famille à un paramètre de géodésiques de type temps 
$D_t\subset \R^{2,1}$, et la famille à un paramètre correspondante $(n(t), v(t))\in T\HH^2$,
pour $t\in [0,1]$. Un calcul facile montre qu'il n'y a pas d'intersection au premier ordre entre 
$D_0$ et $D_t$ (quand $t\to 0$) si et seulement si 
$$ \left\langle n'(0), \frac{\nabla v(t)}{dt}_{|t=0}\right\rangle\neq 0~, $$
où $\langle,\rangle$ désigne le produit scalaire riemannien du plan hyperbolique et
$\nabla$ est sa connexion de Levi-Civita.

Si maintenant $v$ est un champ de vecteurs $\epsilon$-contractant du plan
hyperbolique, alors la dérivée de Lie de la métrique hyperbolique $h$ sous $v$ est
$$ \cL_vh(x,x)=2\langle x,\nabla_xv\rangle $$
si bien que les géodésiques de type temps correspondant à $v$ sont localement
disjointes. 
On peut poursuivre ce raisonnement et montrer que, si $v$ est uniformément contractant,
alors il détermine un feuilletage de $\R^{2,1}$ 
par des géodésiques de type temps.

On peut ensuite procéder comme dans $AdS_3$ pour montrer que l'existence
d'un champ de vecteurs $v$ qui est $(j,\tau)$-équivariant $\epsilon$-contractant assure la
propreté de l'action $(j,\tau)$ de $\Gamma$ sur $\R^{2,1}$, lorsque 
$j$ est convexe co-compacte. En effet on associe à $v$ un feuilletage 
$\cF_0$ de $\R^{2,1}$ par des  géodésiques de type temps, feuilletage qui
est invariant sous l'action de~$\Gamma$ lorsque $v$ est $(j,\tau)$-équivariant.
Puis on montre
que l'action de $\Gamma$ sur l'espace des feuilles est propre. On en déduit
ensuite que l'action de $\Gamma$ sur $\R^{2,1}$ est elle aussi propre.

\section{Déformations en bandelettes et plans croches}
\label{sc:defos-croches}

On va illustrer ici la relation entre déformations en bandelettes (resp.
déformations en bandelettes infinitésimales) et domaines fondamentaux
bordés par des plans croches dans $AdS_3$ (resp. dans $\R^{2,1}$), 
relation qui conduit aux preuves des théorèmes \ref{tm:croches-ads} et 
\ref{tm:croches-minko}. 

\subsection{Bandelettes infinitésimales et plans croches dans Minkowski}

On décrit d'abord le cas plat, qui est à la fois un peu plus simple et 
surtout plus facile à se représenter.

Considérons d'abord le cas le plus simple d'une déformation par bandelette
infinitésimale agissant sur $\HH^2$, sans action de groupe, avec une seule
bandelette infinitésimale. Il s'agit donc
simplement d'un champ de vecteurs $v$ discontinu sur $\HH^2$, qui est~:
\begin{itemize}
\item nul à gauche d'une géodésique orientée $\delta\subset \HH^2$, 
\item égal, à droite de $\delta$, à une translation hyperbolique
infinitésimale $\alpha$ d'axe orthogonal à $\delta$.
\end{itemize}
On note $c$ l'intersection de $\delta$ et de l'axe de $\alpha$, et $u$
le vecteur tangent à $\delta$ en $c$ de longueur égale à la longueur de
$v(c)$, vu comme un vecteur dans 
$\R^{2,1}$. Ainsi, à gauche de $\delta$, $v$~s'écrit simplement sous la
forme $v(x)=u\times x$, où $\times$ est le produit vectoriel naturel de
$\R^{2,1}$.

On a vu au $\S$\ref{ssc:temps-minko} qu'on peut associer à $(x,v(x))\in T\HH^2$
la géodésique de type temps dirigée par $x$ et passant par $x\times v(x)$.
Ici $x\times v(x)=x\times (u\times x)$ est simplement le projeté orthogonal 
de $u$ sur le plan orthogonal à $x$, et on en déduit la description suivante~:
\begin{itemize}
\item lorsque $x$ est à gauche de $\delta$, on associe à $(x,v(x))$ 
la droite passant par $0$ parallèle à $x$, 
\item si $x$ est à droite de $\delta$, on associe à $(x,v(x))$ la droite
passant par $u$ et parallèle à~$x$.
\end{itemize}

On voit ainsi apparaître naturellement deux plans croches gauches, dont les 
centres sont respectivement en $0$ et en $u$ et les parties centrales sont
parallèles à $\delta$. On remarque que ces deux 
plans croches bordent des régions disjointes de $\R^{2,1}$. (Leurs parties
centrales sont contenues dans un même plan de type temps contenant $0$ et $u$.)

Considérons maintenant une surface hyperbolique convexe co-compacte $S$
munie d'une déformation en bandelette infinitésimale. On peut relever 
cette déformation en un champ de vecteurs $v$ sur $\HH^2$, vu comme
le revêtement universel de $S$. Alors $v$ est discontinu le long d'une
famille (infinie) de géodésiques disjointes, et la description ci-dessus s'applique
pour chacune des géodésiques de cette famille. On en déduit une famille de
plans croches gauches --- deux plans croches pour chaque géodésique, chacun étant
l'image de l'autre par une translation --- et on peut vérifier sans 
grande difficulté que des plans croches associés à des géodésiques distinctes
sont disjoints.

Il reste à voir qu'on peut sélectionner une partie de ces plans
croches pour obtenir le bord d'un domaine fondamental pour l'action 
du groupe fondamental de $S$ sur $\R^{2,1}$, il faut pour cela choisir une
famille de géodésiques qui bordent un domaine fondamental de $\HH^2$ pour la
représentation d'holonomie de $S$.

\subsection{Bandelettes et plans croches dans $AdS_3$}
\label{ssc:defos-croches-ads}

On peut reprendre le même schéma pour les déformations non infinitésimales, 
et considérer d'abord une déformation en bandelette très simple de $\HH^2$
avec une seule bandelette. C'est une application discontinue $\phi$ de $\HH^2\setminus
\delta$ dans $\HH^2$, où $\delta$ est encore une géodésique orientée, qui est~:
\begin{itemize}
\item l'identité à gauche de $\delta$,
\item égale, à droite de $\delta$, à une translation hyperbolique $\alpha$ d'axe
orthogonal à $\delta$, qui intersecte $\delta$ en un point $c$.
\end{itemize}
On appelle $B$ la bandelette comprise entre $\delta$ et $\alpha(\delta)$.

Pour tout $x\in \HH^2\setminus B$, on associe à $(x,\phi(x))$ une géodésique
temps $\gamma(x)$ de $AdS_3$ en suivant la construction de la fin du $\S$\ref{ssc:temps-ads},
en identifiant $\HH^2$ au plan $H_0\subset AdS_3$ qui y intervient. On vérifie
sans difficulté, en utilisant l'une des identifications entre géodésiques de type
temps et paires de points de $\HH^2$ du $\S$\ref{ssc:ads}, que~:
\begin{itemize}
\item si $x$ est à gauche de $\delta$, alors $\gamma(x)$ est la géodésique
orthogonale à $H_0$ passant par~$x$. Ces géodésiques feuillettent une moitié du cône de
lumière $C(c_0)$ de $c_0=H_0^*$, le point dual de $H_0$, délimitée par un plan de type temps
$T$ contenant $c_0$ et $\delta$~;
\item si $x$ est à droite de $\delta$, $\gamma(x)$ est une géodésique de type temps
passant par un point~$c_1$ de $T$, obtenu à partir de $c_0$ par une translation 
d'axe orthogonal à la géodésique joignant $c_0$ et $c$ et de longueur égale à la
largeur de la bandelette. Ces géodésiques feuillettent une moitié du cône de lumière
$C(c_1)$ de $c_1$ délimitée par $T$, mais de l'autre coté de la région feuilletée
par les $\gamma(x)$ pour $x$ à gauche de $\delta$.
\end{itemize}
Comme dans le cas Minkowski, on voit donc apparaître deux plans croches gauches~:
celui dont la partie centrale est l'intersection avec $T$ de $C(c_0)$, et celui
dont la partie centrale est l'intersection avec $T$ de $C(c_1)$. Ces deux plans
croches bordent des domaines disjoints de $AdS_3$. 

Considérons maintenant une déformation en bandelettes d'une surface convexe
co-compacte $S$. On peut relever cette déformation en une déformation en 
bandelettes de $\HH^2$, équivariante sous deux représentations. Cette 
déformation se produit le long d'une famille (infinie) de géodésiques de $\HH^2$, 
avec une bandelette introduite le long de chaque géodésique. On associe à
chaque bandelette une paire de plans croches dans $AdS_3$, et deux plans
croches associés à des géodésiques distinctes sont disjoints. En
sélectionnant une partie de ces plans croches comme dans le cas de l'espace de 
Minkowski, on obtient un domaine fondamental
pour l'action sur $AdS_3$ de $(j,\rho)$, où $\rho$ est la représentation
d'holonomie de $S$ et $j$ est la représentation d'holonomie de la structure
hyperbolique obtenue après la déformation en bandelettes.

\bigskip\bigskip\bigskip\bigskip\bigskip\bigskip\bigskip\bigskip\bigskip\bigskip



\newcommand{\etalchar}[1]{$^{#1}$}
\def\cprime{$'$} \def\cprime{$'$}

\end{document}